\documentclass[reqno,10pt,article]{amsart}
\usepackage{amsmath,amsxtra,amssymb,latexsym, amscd,amsthm}
\usepackage[unicode]{hyperref}
\usepackage{array,tabularx,longtable,multicol,indentfirst,fancybox,color}
\usepackage{cases}
\usepackage{graphicx}
\usepackage{tikz}
\usepackage{pgfplots}
\usepackage{multicol}
\usepackage{mathrsfs}
\usepackage{enumerate}

\usetikzlibrary{calc, angles, quotes, intersections, positioning, patterns}

\advance\hoffset-1truecm\relax
\newtheorem{property}{Property}
\newtheorem{theorem}{Theorem}[section]

\newtheorem{lemma}{Lemma}[section]

\newtheorem{remark}{Remark}

\numberwithin{equation}{section}

\def\R{\mathbb{R}}

\begin{document}
	\title[$\mathcal{O}_{\alpha}$-transformation and its uncertainty principles]{$\mathcal{O}_{\alpha}$-transformation and its uncertainty principles}
	\date{\today. Minor Revision: References have been listed in the order they appear in the content and in accordance with the ITSF journal.}		

\author[Lai Tien Minh \& Trinh Tuan]{Lai Tien Minh$^{1}$ and Trinh Tuan$^{2*}$}	
\date{18 February 2026, Accepted by \textsc{Integral Transforms and Special Functions}} 
\thanks{Published online: 04 March 2026. \url{https://doi.org/10.1080/10652469.2026.2635619}}	
\thanks{\textit{$^*$Corresponding author:} \textsc{Trinh Tuan}}

	\maketitle
	\begin{center}
	\scriptsize	
		$^1$Department of Mathematics, Hanoi Architectural University,\\ Hanoi, Vietnam.\\
		E-mail: \texttt{minhlt@hau.edu.vn}\\
		$^2$Department of Mathematics, Faculty of Natural Sciences, Electric Power University,\\ 235-Hoang Quoc Viet Rd., Hanoi, Vietnam.\\
		E-mail: \texttt{tuantrinhpsac@yahoo.com}
	\end{center}			
\begin{abstract}
In this paper, we introduce a family of integral transforms, denoted by \(\mathcal{O}_{\alpha}\), and constructed via kernel fusion of the fractional Fourier transform (FRFT) with angle \(\alpha \notin \pi \mathbb{Z}\). We demonstrate that the \(\mathcal{O}_{\alpha}\)-transformation constitutes a well-defined integral operator by establishing its basic operational properties. Besides, we survey various mathematical aspects of the uncertainty principles for the $\mathcal{O}_{\alpha}$-transform, including Heisenberg's inequality, logarithmic uncertainty inequality, local uncertainty inequality, Hardy's inequality, Pitt's inequality, and Beurling-H{\"o}rmander's theorem.


\end{abstract}	
\section{Introduction}\label{sec1}
\subsection{Background}
\noindent The fractional Fourier transform (abbreviated as  FRFT) is a family of linear transformations generalizing the Fourier transform. It can be regarded as the Fourier transform to the $n$-th power, where $n$ does not need to be an integer, thus, it can transform a function to any intermediate domain between time and frequency. Following \cite{Almeida1994fractional}, the FRFT with angle $\alpha$ is defined in $L^{1}(\mathbb{R})$ with the help of the transformation kernel $K_{\alpha}$ and given by
\begin{equation}\label{eq0}
\mathcal{F}_{\alpha}[f](s)=\int_{\mathbb{R}}K_{\alpha}(t, s)f(t)dt 
\end{equation}
where 
$$
K_{\alpha}(t, s)= \begin{cases}\frac{c(\alpha)}{\sqrt{2 \pi}}  e^{\left\{i a(\alpha)\left(t^{2}+s^{2}-2 b(\alpha)st\right)\right\}}, & \text { if } \alpha \text { is not a multiple of } \pi \\ \delta(t-s), & \text { if } \alpha \text { is a multiple of } 2 \pi \\ \delta(t+s), & \text { if } (\alpha+\pi) \text { is a multiple of } 2 \pi, \end{cases}
$$
with
$
a(\alpha)=\frac{\cot \alpha}{2}, b(\alpha)=\sec \alpha, c(\alpha)=\sqrt{1-i \cot \alpha}.
$ To simplify the notation, throughout this paper, we denote the constants $a(\alpha), b(\alpha)$ and $c(\alpha)$ as $a, b$ and $c$, respectively. 
For $\alpha \in 2 \pi \mathbb{Z}$, the FRFT becomes the identity; for $\alpha+\pi \in 2 \pi \mathbb{Z}$, it is the parity operator. Therefore, from now on we shall confine our attention to $\mathcal{F}_{\alpha}$ for $\alpha \notin \pi \mathbb{Z}$. In the special case angle $\alpha=\pi/2$, we derive the well-known Fourier transform $\mathcal{F}_{\frac{\pi}{2}}:=F$.
In this paper, we introduce  a family of $\mathcal{O}_{\alpha}$-transformation for $L_1$-Lebesgue integrable functions over $\R$ as follows:
\begin{equation}\label{equation1.1}
\mathcal{O}_{\alpha}[f](s):=\int_{\mathbb{R}}\frac{\left[K_{\alpha}(t, s)+ zK_{\alpha}(t,-s)\right]}{2}f(t)dt,
\end{equation}
where $f(t)$ is a function defined on $\mathbb R$. 

\begin{table}[h]
	\centering
	\renewcommand{\arraystretch}{0.8}
	\begin{tabular}{|c|c|c|}
		\hline
		$\alpha$& $z$& $\mathcal{O}_{\alpha}$-transform\\
		\hline
			
		$-\frac{\pi}{2}$&$0$& The Fourier tranform\\
		\hline
				$\frac{\pi}{2}$&$0$& The inverse Fourier Transform\\
		\hline
				$\frac{\pi}{2}$&$1$& The Fourier-cosine transform\\
		\hline
				$\frac{\pi}{2}$&$\pm i$&The Hartley transform\\
		\hline
				$\alpha$&$0$& The fractional Fourier transform\\
         \hline
	\end{tabular}\vskip0.2cm
	\caption{Some of special cases of $\mathcal{O}_{\alpha}$-transform.}\label{table1}
\end{table}
\noindent Some special cases of $\mathcal{O}_{\alpha}$ are enumerated in the Table. \ref{table1}. Note that, if $z=0$, then the transform \eqref{equation1.1} reduces to the FRFT's case.  Therefore, to ensure the newly defined operator is bounded in the $L^2(\mathbb{R})$ space,  in this paper, we only consider the case where $z=\rho i\,,\rho\in\mathbb{R}/\{0\}$.
Then,  the transform \eqref{equation1.1} is not a special case of FRFT \eqref{eq0},  the Dunkl Transform \cite{Rosler1999Bulletin},  as well as the linear canonical transform \cite{canonical}.  Furthermore,  Figure.~\ref{table1} below also displays a comparison between the Fourier, Hartley,  and $\mathcal{O}_{\alpha}$ of the normalized Gaussian function.

\begin{figure*}[h]
	\centering
	\includegraphics[height=9cm,width=14cm]{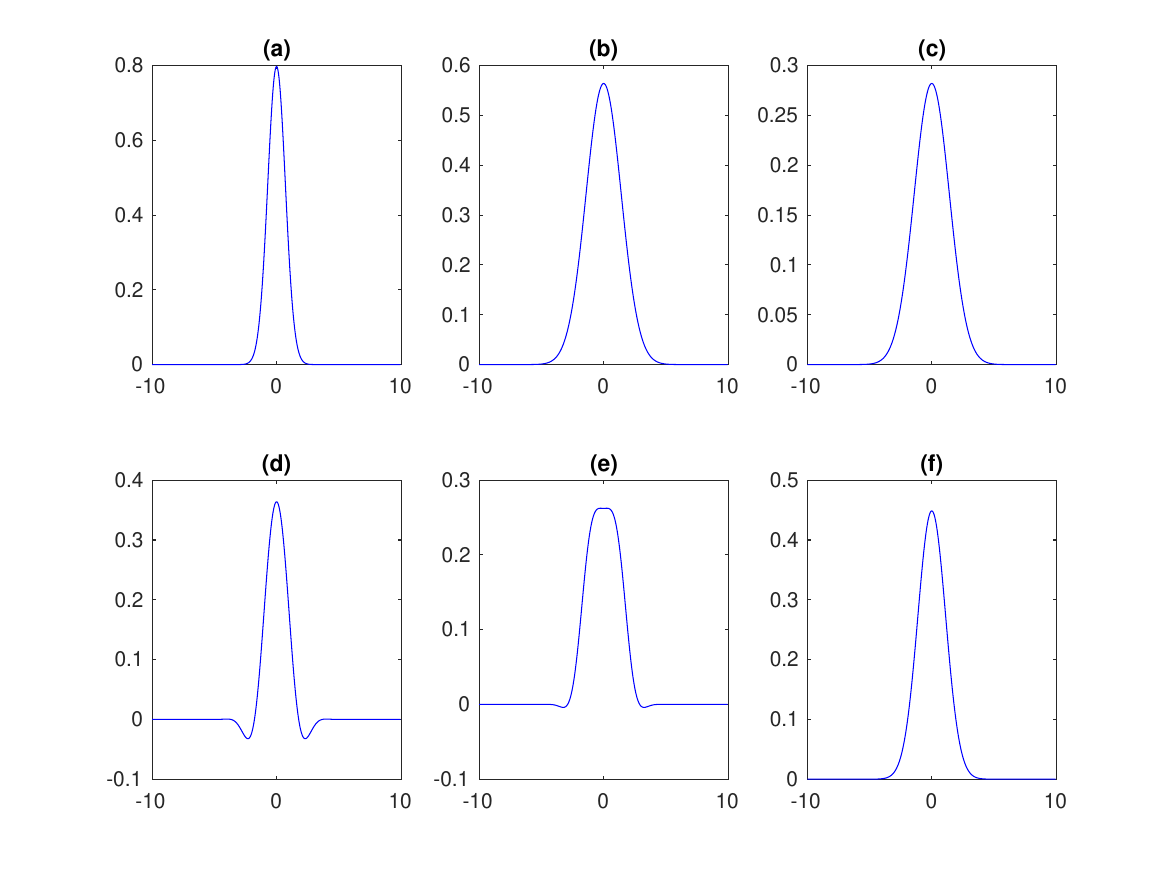}
	\caption{The $\mathcal{O}_{\alpha}$ of the Gaussian function: ({\bf a}) The Gaussian function $f(t)=\sqrt{\frac{2}{\pi}} e^{-t^{2}}$.  ({\bf b}) The Fourier of $f(t)$.  ({\bf c}) The Hartley of $f(t)$.  ({\bf d}) The real part of $\mathcal{O}_{\frac{\pi}{4}}[f](s)$.  ({\bf e}) The imaginary part of $\mathcal{O}_{\frac{\pi}{4}}[f](s)$.  ({\bf f}) The absolute value of $\mathcal{O}_{\frac{\pi}{4}}[f](s)$.}\label{Gauss}
\end{figure*}

\noindent Alternatively,  the transform \eqref{equation1.1} can be explicitly written as
\begin{equation*}
\mathcal{O}_{\alpha}[f](s)=\frac{c(\alpha)}{\sqrt{2 \pi}}e^{ias^{2}}\int_{\mathbb{R}}\frac{\left[e^{-ist\csc\alpha}+ ze^{ist\csc\alpha}\right]}{2}e^{iat^{2}}f(t)dt.
\end{equation*}
\subsection{Organization}

The first purpose of this paper is to establish some basic operator properties for \eqref{equation1.1} to ensure its validity and applicability. Specifically, in Section \ref{sec2}, we prove the $\mathcal{O}_{\alpha}$-transformation maps 
$ L^{1}(\mathbb{R})$ to 
$ C_{0}(\mathbb{R})$ via Riemann-Lebesgue's lemma and  giving Parseval's identity for this transform.  Boundedness of the operator \eqref{equation1.1} has also been pointed out by us through Hausdorff-Young inequality. Besides, we demonstrate a form of the Pitt-type inequality for $\mathcal{O}_{\alpha}$-transform.

The second goal is also the most important of this article. We know that the	 classical uncertainty principle of the Fourier transforms states that a non-zero function and its Fourier
transform cannot both be sharply localized. In other words, it is not possible that the widths
of graphs of $|f (x)|^2$ and $|(Ff )(\xi)|^2$ can both be made arbitrarily \textquotedblleft small\textquotedblright. Concept of smallness has taken different interpretations in different contexts  \cite{cowling1983}, and an outstanding result in \cite{Bonami2003hermite} showed such impossibility when \textquotedblleft smallness\textquotedblright
is interpreted as sharp pointwise or integrable decay. The classical Heisenberg–Weyl uncertainty principle has been extended to other transforms such as Dunkl transform \cite{Rosler1999Bulletin,soltani2013general} and Hankel transform \cite{rosler1999uncertainty,tuan2007uncertainty}. Besides, Heisenberg's uncertainty principle is also key in quantum mechanics \cite{Heisenberg1927,Heisenberg1949physical}, because if we know everything about where a particle is located (the uncertainty of position is small), we know nothing about its momentum (the uncertainty of momentum is large), and vice versa. Its extended versions and variants have important applications in harmonic analysis and physics \cite{Donoho1989uncertainty,folland1997uncertainty,namias1980fractional,Stern2008uncertainty}.  For uncertainty principles for the fractional Fourier transform and other integral operators, while not being exhaustive,  we refer the reader to \cite{Revise3,Revise4,Revise5,Revise2}. Being directly influenced by the above-mentioned, in Section \ref{sec3}, we are motivated to obtain uncertainty principles for the $\mathcal{O}_{\alpha}$-transform, including Heisenberg's inequality, logarithmic uncertainty inequality, and local uncertainty inequality. Hardy's uncertainty inequality and Beurling-H{\"o}rmander type theorem for operator \eqref{equation1.1} are also established at the end of the article.

\section{Operational properties of $\mathcal{O}_{\alpha}$-transform}\label{sec2}
\noindent  This section is devoted to studying some operational properties of the $\mathcal{O}_{\alpha}$-transform \eqref{equation1.1}. For $L^p(\mathbb{R})$, $1\leq p<\infty$, we will use the norm $\|\cdot\|_p$ defined by
$\|f\|_{p}:=\big(\int _{\mathbb{R}}\vert f(y)\vert^pdy\big)^{1/p}.$ The space of all continuous bounded functions over $\mathbb{R}$ that vanish at infinity is denoted by $C_0(\mathbb{R})$. It is a Banach space with respect to the natural supremum-norm $\|\cdot\|_{\infty}$. Property \ref{tc1} follows immediately from Theorems 7.7 in \cite{Rudin1991}.
\begin{property}\label{tc1}
$\mathcal{O}_{\alpha}[f]$ is a bounded linear operator from $L^1(\mathbb{R}) \rightarrow C_0(\mathbb{R})$. In particular, if $f$ be an $L_1$-Lebesgue integrable function over $\mathbb R$, then $\mathcal{O}_{\alpha}[f]$ belongs to $C_0(\mathbb{R})$. Moreover, the following estimation holds

\begin{align}\label{YH1}
\Vert\mathcal{O}_{\alpha}[f]\Vert_{\infty}\leq \frac{\vert\csc\alpha\vert}{\sqrt{2\pi}}\sqrt{\frac{1+\vert z\vert^2}{2}}\cdot\Vert f\Vert_1.
\end{align}
\end{property}
\begin{proof}
A direct computation yields
	$$
	\frac{e^{-ist\csc\alpha}+ ze^{ist\csc\alpha}}{2}=\frac{1}{2}\left[(1+z)\cos\left(st\csc\alpha\right)+ i(z-1)\sin\left(st\csc\alpha\right)\right].
	$$
	Since $\vert (1+z)\cos\left(st\csc\alpha\right)+ i(z-1)\sin\left(st\csc\alpha\right)\vert\leq \sqrt{2(1+|z|^2)},$  we arrive at the following estimation
	$$\begin{aligned}
	\| \mathcal{O}_{\alpha}[f]\ \|_{\infty}&=\sup_{s \in \mathbb{R}}\left\vert(\mathcal{O}_{\alpha}[f](s)\right\vert\\
	&=\sup_{s \in \mathbb{R}}\frac{\vert\csc\alpha\vert}{\sqrt{2\pi}}\left\vert\int_{\mathbb{R}}\frac{\left[e^{-ist\csc\alpha}+ ze^{ist\csc\alpha}\right]}{2}e^{iat^{2}}f(t)dt\right\vert\\
	&\leq \sup_{s \in \mathbb{R}}\frac{\vert\csc\alpha\vert}{\sqrt{2\pi}}\int_{\mathbb{R}}\frac{1}{2} \vert(1+z)\cos\left(st\csc\alpha\right)+ i(z-1)\sin\left(st\csc\alpha\right)\vert \cdot\vert e^{iat^{2}}f(t)\vert dt\\
	&\leq \frac{\vert\csc\alpha\vert}{\sqrt{2\pi}}\sqrt{\frac{1+\vert z\vert^2}{2}}\Vert f\Vert_1.
	\end{aligned}$$
	In addition,  based on the notation $\widetilde{f}(t):=e^{iat^{2}}f(t)$,  we conclude that
	\begin{equation}\label{YH2}
	\mathcal{O}_{\alpha}[f](s)=\frac{c(\alpha)}{2}e^{ias^{2}}F[\widetilde{
		f}](s\csc\alpha)+\frac{c(\alpha)}{2}e^{ias^{2}}F[\widetilde{
		f}](-s\csc\alpha),
	\end{equation}
	where $F$ is the Fourier transform and defined as $F[f](s)= \frac{1}{\sqrt{2\pi}}\int_{\mathbb{R}} f(t)e^{its}dt$.  Since $\vert\widetilde{f}(t)\vert=\vert f(t)\vert$, by combining equation \eqref{YH2} and using Riemann-Lebesgue lemma for  Fourier transform  \cite{Sogge1993fourier},  we infer $\mathcal{O}_{\alpha}[f] \in C_0(\mathbb{R}).$
\end{proof}
To deduce the next property, we need the following useful lemma.
\begin{lemma}[Modified Titchmarsh's representation]\label{tm} Let $f \in L^2(\mathbb{R})$ and assume that $f$ has bounded variation at a point $x$. Then	
	$$
	\frac{1}{2}(f(x+0)+f(x-0))=\lim _{\lambda \rightarrow \infty} \frac{1}{\pi} \int_{\mathbb{R}} f(t) \frac{\sin (\lambda(x-t))}{x-t} d t
	$$\end{lemma}
\begin{proof}	
Since $f \in L^2(\mathbb{R})$ and $(1+|x|)^{-1} \in L^2(\mathbb{R})$, Cauchy-Schwarz implies $\int_{\mathbb{R}} \frac{|f(t)|}{1+|t|} d t$ is finite. Hence $f(t) /(1+|t|) \in L^1(\mathbb{R})$, which is precisely the condition ensuring that Titchmarsh's singular integral representation applies; (see \cite{Titchmarsh1986}, Theorem 12, p. 25). Because $f$ has bounded variation at $x$, both one-sided limits exist, and the claimed identity follows. \end{proof}

\noindent Let $D^k =(d^k \textfractionsolidus dx^k)$ denote the differential operator for $k \in \mathbb{N}$. The Schwartz class is a space of rapidly decreasing functions over $\mathbb{R}$ denoted by $\mathcal{S}(\mathbb{R})=\{f \in C^{\infty}(\mathbb{R}): \|f \|_{k,m} <+\infty \}, \forall k,m \in \mathbb{N}$, where $C^{\infty}(\mathbb{R})$ is a set of smooth functions from $\mathbb{R}\rightarrow \mathbb{C}$ and $\|f \|_{k,m}=\sup_{x\in \mathbb{R}} |x^k D^m f(x)|.$ Schwartz class $\mathcal{S}(\mathbb{R})$ is a dense subset of $L^2 (\mathbb{R})$ \cite{AdamsFournier2003}. 
Denoting the usual inner product in $ L^2(\mathbb{R})$ given by $\langle f, g\rangle=\int_{\mathbb{R}}f(t)\overline{g(t)}dt$,  we deduce Parseval type identity for operator \eqref{equation1.1}  as follows:

\begin{property}\label{tc2}
For any two function $f, g \in L^2 (\mathbb R)$ one has
	\begin{equation}\label{parseval1}
	\langle \mathcal{O}_{\alpha}f,\mathcal{O}_{\alpha}g\rangle =\frac{1}{4}(1+\vert z\vert^2)\langle f,g\rangle.
	\end{equation}
	In the special case of $f=g$ , we then have
	\begin{equation}\label{parseval2}
	\| \mathcal{O}_{\alpha}f\|_2^2=\frac{1}{4}(1+\vert z\vert^2)\|f\|_2^2.
	\end{equation}
\end{property}
\begin{proof}
	As far as we know $\mathcal{S}(\mathbb{R})$ is dense in $L_2(\mathbb{R})$,  we then only need to prove \eqref{parseval1} hold true for $f, g \in \mathcal{S}(\mathbb{R})$.  We realize that $z+\overline{z}=0$.  Hence, with the help of Def. \eqref{equation1.1}, we then have	
	\begin{equation}\label{H3}
\begin{aligned} &\langle \mathcal{O}_{\alpha}f,\mathcal{O}_{\alpha}g\rangle = \int_{\mathbb{R}} \mathcal{O}_{\alpha}[f](s) \overline{\mathcal{O}_{\alpha}[g](s)} ds\\ &=\frac{\csc\alpha}{8\pi}\int_{\mathbb{R}^3}\left[e^{-ist\csc\alpha}+ ze^{ist\csc\alpha}\right]\left[e^{isx\csc\alpha}+ \overline{z}e^{-isx\csc\alpha}\right] \widetilde{f}(t)\overline{\widetilde{g}(x)} dt dx ds\\ &=\frac{\csc\alpha}{8\pi}\int_{\mathbb{R}^2} \left\{\int_{\mathbb{R}}\left[e^{is(x-t)\csc\alpha} +\overline{z} e^{-is(x+t)\csc\alpha} +ze^{is(x+t)\csc\alpha} +\vert z\vert^2 e^{is(t-x)\csc\alpha}\right] ds\right\}\widetilde{f}(t)\overline{\widetilde{g}(x)}dt dx\\ &=\frac{\csc\alpha}{8\pi}\int_{\mathbb{R}^2}\widetilde{f}(t)\overline{\widetilde{g}(x)}\times \bigg\{\int_{\mathbb{R}}\bigg((1+\vert z\vert^2)\cos[s(x-t)\csc\alpha]+i(z-\overline{z})\sin[s(x+t)\csc\alpha]\\ &+i(1-\vert z\vert^2)\sin[s(x-t)\csc\alpha]+(z+\overline{z})\cos[s(x+t)\csc\alpha]\bigg)ds   \bigg\} dt dx\\ &=\frac{\csc\alpha}{8\pi}\lim_{\lambda\to\infty}\!\int_{\mathbb{R}^2}\widetilde{f}(t)\overline{\widetilde{g}(x)}\times \bigg\{\int_{-\lambda}^\lambda \bigg((1+\vert z\vert^2)\cos[s(x-t)\csc\alpha]+i(z-\overline{z})\sin[s(x+t)\csc\alpha]\\ &+i(1-\vert z\vert^2)\sin[s(x-t)\csc\alpha]\bigg)ds\bigg\} dt dx. \end{aligned}
	\end{equation}
	Thanks to formula $$\displaystyle\int_{-\lambda}^{\lambda} \big(\sin[s(x+t)\csc\alpha]\big)ds=\displaystyle\int_{-\lambda}^{\lambda} \big(\sin[s(x-t)\csc\alpha]\big)ds=0,$$ then  \eqref{H3} becomes
	$$
	\langle \mathcal{O}_{\alpha}f,\mathcal{O}_{\alpha}g\rangle =\frac{1}{4}(1+\vert z\vert^2)\int_{\mathbb{R}}\widetilde{f}(t)\left\{\lim_{\lambda\to\infty} \frac{1}{\pi}\int_{\mathbb{R}} \overline{\widetilde{g}(x)} \frac{\sin \lambda(x-t)\csc\alpha}{x-t} dx \right\}  dt.
	$$
	\noindent	By virtue of Lemma \eqref{tm},  the relation above turns out to be
	
	\begin{align*}
	\langle \mathcal{O}_{\alpha}f,\mathcal{O}_{\alpha}g\rangle&=\frac{1}{4}(1+\vert z\vert^2)\int_{\mathbb{R}} \widetilde{f}(t)\overline{\widetilde{g}(x)}dt \\&=\frac{1}{4}(1+\vert z\vert^2)\int_{\mathbb{R}} f(t)\overline{g(t)} dt =\frac{1}{4}(1+\vert z\vert^2)\langle f,g\rangle.
	\end{align*}

\noindent	The proof is complete.
\end{proof}
	\noindent In what follows, we will show the boundedness for the operator \eqref{equation1.1} on the dual space.  The techniques used in the proof of our theorem come from \cite{Duoandikoetxea01,MinhTuan25}. 
\begin{property}[Hausdorff-Young type inequality for $\mathcal{O}_{\alpha}$-transform]
Let $p\in[1,2]$ and  $p_1$ denote the conjugate exponent of $p$. Then, we have the following estimation 	
\begin{equation}\label{eq2.6}
\|\mathcal{O}_{\alpha}f\|_{p_1} \leq \left(\frac{|\csc\alpha|}{2\sqrt{\pi}}\sqrt{\frac{1+\vert z\vert^2}{2}}\right)^{\zeta} \|f\|_{p},\quad \forall f\in L_p(\mathbb{R}),
\end{equation}
where,  $\zeta=\frac{2}{p}-1$.	
\end{property}
\begin{proof}
	Since $f$ belongs to $L_p(\mathbb{R})$,  it is straightforward to get $\widetilde{f}\in L_p(\mathbb{R})$, $p\in[1,2]$.  Moreover, making use of Babenko–Beckner's inequality \cite{Beckner1975},  we obtain  
	$F[\widetilde{
		f}](s\csc\alpha) \in L_{p_1}(\mathbb{R})$ and  $F[\widetilde{
		f}](-s\csc\alpha) \in L_{p_1}(\mathbb{R}), $
	where $\frac{1}{p}+\frac{1}{p_1}=1,$ with $p \in [1,2]$. It is easy to verify that $|e^{iat^{2}}| =1$.  Now,  owing to \eqref{YH1}, we derive $\mathcal{O}_{\alpha}[f]\in L_{p_1}(\mathbb{R}).$
	On the other hand,  Property \ref{tc1} gives us that the  $\mathcal{O}_{\alpha}$-transform is a bounded linear operator from $L_1(\mathbb{R})\to L_\infty(\mathbb{R})$.
	Moreover,  with the aid of \eqref{parseval2},  we deduce $\mathcal{O}_{\alpha}$ is a bounded from $L_2(\mathbb{R} )\to L_2(\mathbb{R})$.
	Applying the Riesz–Thorin interpolation theorem \cite{Duoandikoetxea01},  we obtain $\mathcal{O}_{\alpha}$ is a bounded linear operator from $L_p(\mathbb{R})\to L_{p_1}(\mathbb{R})$, where $p\in[1,2]$, $\frac{1}{p}+\frac{1}{p_1}=1$ and $\frac{1}{p}=\frac{\zeta}{1}+\frac{1-\zeta}{2}$ (and then $\zeta=\frac{2}{p}-1$). Therefore,  the desired estimation can be achieved as follows
	$$
	\|\mathcal{O}_{\alpha}f\|_{p_1} \leq \left(\frac{|\csc\alpha|}{2\sqrt{\pi}}\sqrt{\frac{1+\vert z\vert^2}{2}}\right)^{\zeta} \|f\|_{p}.
	$$
	This indicates that \eqref{eq2.6} is proved.
\end{proof}

 Next, we recall the classical Pitt's inequality for Fourier transform.
For any  $f \in \mathcal{S}(\mathbb{R}) \subseteq L^{2}(\mathbb{R})$, Pitt's inequality is given as
\begin{equation}\label{Pitt0}
\int_{\mathbb{R}}|s|^{-\lambda}\left|F[f](s)\right|^{2} ds \leq C_{\lambda} \int_{\mathbb{R}}\left|t\right|^{\lambda}\left| f(t)\right|^{2}d t,\ \text{with}\ 0 \leq \lambda<1,  
\end{equation}
where the upper bound constant $C_{\lambda}$ is defined by
\begin{equation}\label{Pitt1}
C_{\lambda}=\pi^{\lambda}\left[\frac{\Gamma\left(\frac{1-\lambda}{4}\right)}{\Gamma\left(\frac{1+\lambda}{4}\right)} \right]^{2}.
\end{equation}

\noindent Beckner's result in \cite{Beckner1995pitt} demonstrates that a sharp form of Pitt’s inequality directly yields a concise proof of the logarithmic uncertainty principle. Moreover, the classical Heisenberg–Weyl inequality can be derived as a corollary of this logarithmic estimate in the setting of the real number field. In \cite{Beckner1995pitt}, employing techniques of rearrangement and symmetrization, Beckner established the sharp form of Pitt’s inequality by utilizing the optimal $L^1$-Young’s inequality for convolution on $\mathbb{R}_+$. Directly motivated by Beckner's works \cite{Beckner1995pitt,Beckner2008pitt}, we present an alternative proof of Pitt’s inequality in the context of the $\mathcal{O}_{\alpha}$-transformation over the entire real line $\mathbb{R}$ as follows.
\begin{property}[Pitt's inequality for $\mathcal{O}_{\alpha}$-transform]\label{tc4}
	For any function $f \in\mathcal{S}(\mathbb{R})$, we have 
	\begin{align}\label{Pitt2}
	|\csc\alpha|^{-\lambda}& \int_{\mathbb{R}}|s|^{-\lambda}\left|\mathcal{O}_{\alpha}[f](s)\right|^{2} ds\leq \frac{1}{4} C_{\lambda}(1+\vert z\vert^2) \int_{\mathbb{R}}|t|^{\lambda}|f(t)|^{2}dt,
	\end{align}
	where $C_{\lambda}$ is given by \eqref{Pitt1}.
\end{property}
\begin{proof}
	Using inequality \eqref{Pitt0} for the function $\widetilde{f}(t)$, we obtain
	$$
	\int_{\mathbb{R}}|s|^{-\lambda}|F[\widetilde{f}](s)|^{2}ds \leq C_{\lambda}\int_{\mathbb{R}}|t|^{\lambda}|\widetilde{f}(t)|^{2}dt.
	$$
	Equivalently
	\begin{equation*}
	|\csc\alpha| \int_{\mathbb{R}}|s\csc\alpha|^{-\lambda}|F[\widetilde{f}](s\csc\alpha)|^{2}ds \leq C_{\lambda} \int_{\mathbb{R}}|t|^{\lambda}|\widetilde{f}(t)|^{2} dt.
	\end{equation*}
	Notice that $\left|\mathcal{F}_{\alpha}[f](s)\right|=\sqrt{|\csc\alpha|} \left|F[\widetilde{f}](s\csc\alpha)\right|$ and $|\widetilde{f}(t)|=|f(t)|$, the inequality \label{Pitt3} yields
	\begin{equation}\label{a}
	|\csc\alpha|^{-\lambda} \int_{\mathbb{R}}|s|^{-\lambda}\left|\mathcal{F}_{\alpha}[f](s)\right|^{2} ds \leq C_{\lambda} \int_{\mathbb{R}}|t|^{\lambda}|f(t)|^{2}dt.
	\end{equation}
	Similarly, we infer 
	\begin{equation}\label{b}
	|\csc\alpha|^{-\lambda} \int_{\mathbb{R}}|s|^{-\lambda}\left|\mathcal{F}_{\alpha}[f](-s)\right|^{2} ds \leq C_{\lambda} \int_{\mathbb{R}}|t|^{\lambda}|f(t)|^{2}dt .
	\end{equation}
	By simple computations, we have $|x+zy|^2=|x|^2+|z|^2\cdot|y|^2+\overline{z}x\overline{y}+z\overline{x}y.$
	Based on the inequality $xy\le\frac{|x|^2+|y|^2}{2}, $ we obtain 
	$$
	|x+zy|^2\leq \left(1+\frac{z+\overline{z}}{2}\right)|x|^2+\left(|z|^2+\frac{z+\overline{z}}{2}\right)\cdot|y|^2.
	$$
	Since $z+\overline{z}=0$, then the estimate above turns to 
	\begin{align}\label{complex}
	|x+zy|^2\leq |x|^2+\vert z\vert^2|y|^2.
	\end{align}
	By combining relation $\mathcal{O}_{\alpha}[f](s)=\frac{1}{2}\left[\mathcal{F}_{\alpha}[f](s)+z\mathcal{F}_{\alpha}[f](-s)\right]$,   inequality \eqref{complex},  equation \eqref{a} and \eqref{b} we obtain 
	\begin{equation}\label{c}
	\begin{aligned}
	|\csc\alpha|^{-\lambda} \int_{\mathbb{R}}|s|^{-\lambda}\left|\mathcal{O}_{\alpha}[f](s)\right|^{2} ds
	&\leq \frac{1}{4}|\csc\alpha|^{-\lambda} \int_{\mathbb{R}}|s|^{-\lambda}\left\vert\mathcal{F}_{\alpha}[f](s)+z\mathcal{F}_{\alpha}[f](-s)\right\vert^2 ds\\
	&\leq \frac{1}{4} C_{\lambda}(1+\vert z\vert^2)\int_{\mathbb{R}}|t|^{\lambda}|f(t)|^{2}dt.
	\end{aligned}
	\end{equation}	
\noindent Inequality \eqref{c} can be called Pitt's inequality for $\mathcal{O}_{\alpha}$-transform.
\end{proof}

\section{Uncertainty principles associated with   $\mathcal{O}_{\alpha}$-transform}\label{sec3}
\subsection{Heisenberg's uncertainty principle for $\mathcal{O}_{\alpha}$-transform}
\begin{theorem}\label{theorem2} 
	If $f\in L^2 (\mathbb{R})$, then 
	\begin{equation}\label{heisenberg2}
	\Vert tf(t)\Vert_2\cdot \Vert s\mathcal{O}_{\alpha}[f](s)\Vert_2\ge |\sin\alpha|\sqrt{\frac{1+\vert z\vert^2}{2}}\Vert f\Vert_2^2.
	\end{equation}
	The equal sign in inequality \eqref{heisenberg2} occurs if and only if
	 $f(u)=Ae^{-\beta u^2}e^{-iu^{2}\cot\alpha},$ where $A,\beta$ are some constants and $\beta>0$.
\end{theorem}
\begin{proof}
	It suffices to prove that $f\in\mathcal{S}(\mathbb{R})$. Since $\mathcal{S}(\mathbb{R})$ is dense in $L^2 (\mathbb{R})$ it suffices to verify the estimate for
	$f \in \mathcal{S}(\mathbb{R})$.  Integrating by parts and using $f, f'\in\mathcal{S}(\mathbb{R})$, we obtain
	\begin{equation*}
	\Vert f \Vert_2^2=-\int_{\mathbb{R} }u\frac{d}{du}|f(u)|^2
	=-\langle uf'(u), f(u)\rangle-\langle uf(u), f'(u)\rangle.
	\end{equation*}
	Based on  the Cauchy-Schwarz inequality and the relation $$|\langle uf'(u), f(u)\rangle|=|{\langle}uf(u), f'(u)\rangle|,$$ we then have
	\begin{equation*}
	\Vert f \Vert_2^2\le 2\int_{\mathbb{R} }|u|\cdot|f(u)|\cdot|f'(u)|du
	\le 2\Vert uf(u)\Vert_2\cdot\Vert f' \Vert_2.
	\end{equation*}
	Thus
	\begin{equation}\label{C-S1}
	\Vert f \Vert_2^2\le 2\Vert uf(u)\Vert_2\cdot\Vert f' \Vert_2.
	\end{equation}
	The equal sign in \eqref{C-S1} occurs if and only if $f(u)=Ae^{-\beta u^2},$ where $\beta$ is a positive constant for $f(u)$ to be in $\mathcal{S}(\mathbb{R})$ \cite{folland1997uncertainty}, this implies that $f'(u)=\beta uf(u)$ for some constants $\beta$. 
	Conversely, if $f(u)=Ae^{-\beta u^2}$ and  $\beta$ is a positive constant for $f(u)$ to be in $\mathcal{S}(\mathbb{R})$ \cite{folland1997uncertainty}.
	Since $\mathcal{S}(\mathbb{R})$ is dense in $L^2(\mathbb{R})$, it suffices to prove this theorem for $f\in\mathcal{S}(\mathbb{R})$. Integrating by parts and using $f, f'\in\mathcal{S}(\mathbb{R})$, we obtain
	\begin{equation*}
	\Vert f \Vert_2^2=-\int_{\mathbb{R} }u\frac{d}{du}|f(u)|^2du.
	\end{equation*}
	Using the relation $|\langle uf'(u), f(u)\rangle|=|{\langle}uf(u), f'(u)\rangle|$,\,
	and the Cauchy-Schwarz inequality, results in 
	\begin{equation}\label{C-S}
	\Vert f \Vert_2^2\le 2\int_{\mathbb{R} }|u|\cdot|f(u)|\cdot|f'(u)|du
	\le 2\Vert uf(u)\Vert_2\cdot\Vert f' \Vert_2.
	\end{equation}
	To begin with, we put
	\begin{align*}
	&\mathcal{H}_{\alpha}[f](s) : = \frac{c(\alpha)}{\sqrt{2\pi}}\int_{\mathbb{R}}\left(e^{-ist\csc\alpha}+ ze^{ist\csc\alpha}\right)f(t)dt. 
	\end{align*}
	By using a similar idea from the proof of the Parseval identity \ref{parseval1},  we can demonstrate that $\mathcal{H}_{\alpha}$ is a continuous operator from $\mathcal{S}(\mathbb{R})$ onto $\mathcal{S}(\mathbb{R})$. In addition, if $f(t), \mathcal{H}_{\alpha}[f](s)  \in \mathcal{S}(\mathbb{R}),$  we deduce that 
\begin{equation}\label{equation2.1}
\begin{aligned}
	&\frac{c(-\alpha)}{\sqrt{2\pi}(1+\vert z\vert^2)}\int_{\mathbb{R}}\left[e^{ist\csc\alpha}+\overline{z}e^{-ist\csc\alpha}\right]\mathcal{H}_{\alpha}[f(v)](t) dt\\
	&=\frac{\csc\alpha}{2\pi(1+\vert z\vert^2)}\int_{\mathbb{R}}\int_{\mathbb{R}}\left[e^{ist\csc\alpha}+\overline{z}e^{-ist\csc\alpha}\right]\left[e^{-itv\csc\alpha}+ ze^{itv\csc\alpha}\right]f(v)dvdt\\
	&=\lim_{\lambda\to\infty}\frac{\csc\alpha}{2\pi(1+\vert z\vert^2)}\int_{\mathbb R}f(v)dv\int_{-\lambda}^{+\lambda}\left[e^{ist\csc\alpha}+\overline{z}e^{-ist\csc\alpha}\right]\left[e^{-itv\csc\alpha}+ ze^{itv\csc\alpha}\right]dt\\
	&=\lim_{\lambda\to\infty}\frac{\csc\alpha}{2\pi}\int_{\mathbb R}f(v)dv\int_{-\lambda}^{+\lambda}\cos \left[ t(s-v)\csc\alpha\right]dt\\
	&=\lim_{\lambda\to\infty}\frac{1}{\pi}\int_{\mathbb R}f(v)\frac{\sin\lambda(s-v)\csc\alpha}{s-v}dv=f(s).
	\end{aligned}\end{equation}
	In other words 
	$$\mathcal{H}_{\alpha}^{-1}\{f(t)\}(s) : = \frac{c(-\alpha)}{\sqrt{2\pi}(1+\vert z\vert^2)}\displaystyle\int_{\mathbb{R}}\big(e^{jst\csc\alpha}- ze^{-ist\csc\alpha}\big)f(t)dt.$$
	Hence, $\mathcal{H}_{\alpha}$ is a linear, continuous and one-to-one mapping from $\mathcal{S}(\mathbb{R})$ onto $\mathcal{S}(\mathbb{R})$ whose inverse is also continuous.  Under the condition $f\in L^1({\mathbb{R}} )$, and let $g\in\mathcal{S}(\mathbb{R}),$ we validate the results $$\int_{\mathbb{R}}f(u)\mathcal{H}_{\alpha}[g](u)du=\int_{\mathbb{R}}g(t)\mathcal{H}_{\alpha}[f](t)dt.$$ On the other hand,  with the aid of above equality and relation \eqref{equation2.1},  we infer that

	\begin{align*}
	\int_{{\mathbb{R}}}\!f(s)\mathcal{H}_{\alpha}[g](s)ds
	&=\frac{c(-\alpha)}{\sqrt{2\pi}(1+\vert z\vert^2)}\int_{{\mathbb{R}}}\!\Bigg[\int_{{\mathbb{R}}}\left(e^{ist\csc\alpha}+\overline{z}e^{-ist\csc\alpha}\right)\mathcal{H}\{g(t)\}(u)du\Bigg]\mathcal{H}_{\alpha}[f](t)dt\\
	&=\int_{{\mathbb{R}}}\!(\mathcal{H}_{\alpha}g)(s)\Bigg[\frac{c(-\alpha)}{\sqrt{2\pi}(1+\vert z\vert^2)}\!\int_{{\mathbb{R}}}\left(e^{ist\csc\alpha}+\overline{z}e^{-ist\csc\alpha}\right)\mathcal{H}_{\alpha}[f](s)dt\Bigg]ds\\
	&=\int_{{\mathbb{R}}}\!\!f_0(s)\mathcal{H}_{\alpha}[f](s)ds,
	\end{align*}
	where $f_0(s)$ defined by $$f_0(s)=\frac{c(-\alpha)}{\sqrt{2\pi}(1+\vert z\vert^2)}\int_{\mathbb{R}}\left(e^{ist\csc\alpha}+\overline{z}e^{-ist\csc\alpha}\right)\mathcal{H}_{\alpha}[f](t)dt.$$
	As mentioned above,  the function $\mathcal{H}_{hg}$ covers all of $\mathcal{S}(\mathbb{R})$ when $g$ runs over  $\mathcal{S}(\mathbb{R})$. Therefore,
	$$\displaystyle\int_{\mathbb{R}}[f_0(s)-f(s)]\Phi(s)ds=0,\ \forall \Phi\in \mathcal{S}(\mathbb{R}).$$  Having in mind that $\mathcal{S}(\mathbb{R})$ is dense in $L^1(\mathbb{R} )$, and $f_0(s)-f(s)=0$  almost every for $s\in{\mathbb{R}}$, we deduce that $f\in L^1(\mathbb{R}),$ and $\mathcal{H}_{\alpha}[f]\in L^1(\mathbb{R}),$ then the following relation 
	\begin{equation*}
	f(s)=\frac{c(-\alpha)}{\sqrt{2\pi}(1+\vert z\vert^2)}\int_{\mathbb{R}}\left(e^{ist\csc\alpha}+\overline{z}e^{-ist\csc\alpha}\right)\mathcal{H}_{\alpha}[f](t)dt
	\end{equation*}
	holds for almost every $s\in{\mathbb{R}}$.
	In the sequel, we set 
	\[\mathcal{H}_{\alpha}^{-1}[g(t)](s):=\frac{c(-\alpha)}{\sqrt{2\pi}(1+\vert z\vert^2)}\int_{\mathbb{R}}\left(e^{ist\csc\alpha}+\overline{z}e^{-ist\csc\alpha}\right)g(t)dt\]
	as the inverse transform of $\mathcal{H}_{\alpha}.$ 
	Let us consider  the following transform
	\begin{align*}
	&\mathscr{Q}\{f(t)\}(s) : =\frac{c(-\alpha)\csc\alpha}{\sqrt{2\pi}(1+\vert z\vert^2)}\int_{\mathbb{R}}\left(e^{ist\csc\alpha}-\overline{z}e^{-ist\csc\alpha}\right)f(t)dt.
	\end{align*}
	Using $\mathcal{H}_{\alpha}^{-1}$ and $\mathcal{H}_{\alpha}$, $f$ can take the form
	$f=\mathcal{H}_{\alpha}^{-1}\{\mathcal{H}_{\alpha}[f]\}.$
	This in turn implies that
	$$
	f'(s)=\mathscr {Q}\big\{t\mathcal{H}_{\alpha}[f](t)\big\}(s).
	$$
	Applying the definition of  inner product on $L^2(\mathbb{R})$, we obtain 
		\begin{align*}
		&\Vert f' \Vert_2^2=\int_{\mathbb{R}}f'(s)\overline{f'(s)}ds=\lim_{\lambda\to\infty}\int_{-\lambda}^{+\lambda}f'(s)\overline{f'(s)}ds\\
		&=\lim_{\lambda\to\infty}\int_{-\lambda}^{+\lambda}\mathscr {Q}\big\{t\mathcal{H}_{\alpha}[f](t)\big\}(s)\cdot\overline{\mathscr {Q}\big\{t\mathcal{H}_{\alpha}[f](t)\big\}(s)}ds\\
		&=\lim_{\lambda\to\infty}\frac{\csc^3\alpha}{4\pi(1+\vert z\vert^2)^2}\int_{-\lambda}^{+\lambda}\int_{{\mathbb R}}t\mathcal{H}_{\alpha}[f](t)dt\int_{{\mathbb R}}\overline{v\mathcal{H}_{\alpha}[f](v)}dv\left[e^{ist\csc\alpha}-\overline{z}e^{-ist\csc\alpha}\right]\left[e^{-isv\csc\alpha}-ze^{isv\csc\alpha}\right]ds\\
		&=\lim_{\lambda\to\infty}\frac{\csc^3\alpha}{4\pi(1+\vert z\vert^2)^2}\int_{{\mathbb R}}t\mathcal{H}_{\alpha}[f](t)dt\int_{{\mathbb R}}\overline{v\mathcal{H}_{\alpha}[f](v)}dv\int_{-\lambda}^{+\lambda}\left[e^{ist\csc\alpha}-\overline{z}e^{-ist\csc\alpha}\right]\left[e^{-isv\csc\alpha}-ze^{isv\csc\alpha}\right]ds\\
		&=\lim_{\lambda\to\infty}\frac{\csc^3\alpha}{4\pi(1+\vert z\vert^2)}\int_{{\mathbb R}}t\mathcal{H}_{\alpha}[f](t)dt\int_{{\mathbb R}}\overline{v\mathcal{H}_{\alpha}[f](v)}dv\int_{-\lambda}^{+\lambda}\cos\left[s(t-v)\csc\alpha\right]ds\\
		&=\lim_{\lambda\to\infty}\frac{\csc^2\alpha}{2\pi(1+\vert z\vert^2)}\int_{\mathbb R}t\mathcal{H}_{\alpha}[f](t)dt\int_{\mathbb R}\overline{v\mathcal{H}_{\alpha}[f](v)}\frac{\sin \left[(t-v)\lambda\csc\alpha\right]}{t-v}dv.
		\end{align*}
	Lemma \ref{tm} yields the following relation 	
		\begin{align*}
		\Vert f' \Vert_2^2&=\frac{\csc^2\alpha}{2(1+\vert z\vert^2)}\int_{{\mathbb R}}t\mathcal{H}_{\alpha}[f](t)\left\{\lim_{\lambda\to\infty}\frac{1}{\pi}\int_{\mathbb R}\overline{v\mathcal{H}_{\alpha}[f](v)}\frac{\sin \left[(t-v)\lambda\csc\alpha\right]}{t-v}dv\right\}dt\\
		&=\frac{\csc^2\alpha}{2(1+\vert z\vert^2)}\int_{\mathbb R}t\mathcal{H}[f](t)\overline{t\mathcal{H}_{\alpha}[f](t)}dt=\frac{\csc^2\alpha}{2(1+\vert z\vert^2)}\Vert t\mathcal{H}_{\alpha}[f](t) \Vert_2^2.
		\end{align*}
	Therefore
	\begin{equation}\label{eqM}
	\Vert f' \Vert_2^2=\frac{\csc^2\alpha}{2(1+\vert z\vert^2)}\Vert t\mathcal{H}_{\alpha}[f](t) \Vert_2^2.
	\end{equation}
	By inserting the relation \eqref{eqM} into inequality \eqref{C-S}, it turn out to be
	\begin{equation}\label{heisenbergH}
	\Vert uf(u)\Vert_2\cdot \Vert s\mathcal{H}_{\alpha}[f](s)\Vert_2\ge |\sin\alpha|\sqrt{\frac{1+\vert z\vert^2}{2}}\Vert f\Vert_2^2.
	\end{equation}
	Obviously
	$$\mathcal{H}_{\alpha}[f(t)e^{it^2\frac{\cot\alpha}{2}}](s)=e^{-is^2\frac{\cot\alpha}{2}}\mathcal{O}_{\alpha}\{f(t)\}(s)\quad 
	\text{and} \quad 
	\big|e^{-\frac{jd}{2b}s^{2}}\big|=\big|e^{\frac{ja}{2b}t^{2}}\big|=1.$$
	Therefore, applying inequality \eqref{heisenbergH} in the case $f(t)e^{it^{2}\frac{\cot\alpha}{2}}$, we deduce \eqref{heisenberg2}. Similarly, the equality occurs in each of \eqref{heisenberg2} if and only if $f(u)=Ae^{-\beta u^2}e^{-iu^2\cot\alpha},$ where $A$ and $\beta$ are constants.  The proof is complete.
\end{proof}
\begin{remark}
\textup{For the case angle $\alpha=\pi/2$ Theorem \ref{theorem2} reduces to the uncertainty inequality pertaining to Fourier transform:	
$$\Vert uf(u)\Vert_2 . \Vert sF\{f(t)\}(s)\Vert_2\ge \frac{1}{\sqrt{2}}\Vert f\Vert_2^2.$$
Considering the normalized Gaussian function $f(t)=\sqrt{\frac{2}{\pi}} e^{-t^{2}}$,  we can compute via the standard-Gaussian integral as
	\begin{equation*}
	\Vert tf(t)\Vert_2^{2}=\int_{\mathbb{R}}|t|^{2}|f(t)|^{2} d t=\frac{2}{\pi} \int_{\mathbb{R}} t^{2} e^{-2 t^{2}} d t=\frac{4}{\pi} \int_{0}^{\infty} t^{2} e^{-2 t^{2}} dt=\frac{1}{2 \sqrt{2 \pi}}. 
	\end{equation*}	
	Moreover
		\begin{equation*}
	\Vert f \Vert_2^{2}=\int_{\mathbb{R}}|f(t)|^{2} d t=\frac{2}{\pi} \int_{\mathbb{R}} e^{-2 t^{2}} d t=\sqrt{\frac{2}{ \pi}}. 
	\end{equation*}	
	Using \eqref{heisenberg2},  the following inequality must be satisfied
	\begin{equation*}\label{eq2.16}
	\Vert s\mathcal{O}_{\alpha}[f](s)\Vert_2^2 \geq 2\sqrt{\frac{2}{\pi}}\sin^2\alpha. 
	\end{equation*} }
\end{remark}

\subsection{Logarithm uncertainty principles for $\mathcal{O}_{\alpha}$-transform}

Beckner's result in \cite{Beckner1995pitt} is a novel addition to the class of quantitative uncertainty inequalities governing the localization of a function $f$ and the $F[f]$ by virtue of a logarithmic estimate obtained from the classical Pitt's inequality. For any $f \in \mathcal{S}(\mathbb{R}) \subseteq L^{2}(\mathbb{R})$,  the Beckner's inequality \cite{Beckner1995pitt}  states that
\begin{equation}\label{Beck0}
\int_{\mathbb{R}} \ln |t||f(t)|^{2}dt+\int_{\mathbb{R}} \ln |s||F[f](s)|^{2}ds \geq\left(\frac{\Gamma^{\prime}(1 / 4)}{\Gamma(1 / 4)}-\ln \pi\right) \int_{\mathbb{R}}|f(t)|^{2} dt
\end{equation}
where $\Gamma$ is the well-known Euler Gamma function.
The inequality \eqref{Beck0} shares a bond with Heisenberg's uncertainty principle, so it is often labeled as the logarithmic uncertainty principle. We point out a version of logarithm uncertainty for $\mathcal{O}_{\alpha}$-transform as follows

\begin{theorem}\label{logarithmic}
	For any $f\in\mathcal{S}(\mathbb{R})$, we have
	\begin{align*}
	 &\frac{1}{4}(1+\vert z\vert^2)\int_{\mathbb{R}} \ln |t||f(t)|^{2}dt+\int_{\mathbb{R}} \ln |s|\left|\mathcal{O}_{\alpha}[f](s)\right|^{2}ds \\&\geq \frac{1}{4}(1+\vert z\vert^2)\left(\frac{\Gamma^{\prime}(1 / 4)}{\Gamma(1 / 4)}+\ln \left\vert\frac{\sin\alpha}{\pi}\right\vert\right) \int_{\mathbb{R}}|f(t)|^{2}dt .
	\end{align*}
\end{theorem}
\begin{proof}
	We have $F[\widetilde{f}](s)=F(s\csc\alpha)$.  
	Also, for every $0 \leq \lambda<1$, we define
	\begin{equation}\label{d}
	M(\lambda)=|\csc\alpha|^{-\lambda} \int_{\mathbb{R}}|s|^{-\lambda}\left|\mathcal{O}_{\alpha}[f](s)\right|^{2} ds-\frac{1}{4}(1+\vert z\vert^2)C_{\lambda} \int_{\mathbb{R}}|t|^{\lambda}|f(t)|^{2}dt .
	\end{equation}
	As can be seen $M(\lambda)\leq 0,\forall \lambda \in[0,1)$. Moreover, using the relation \eqref{parseval1}, we then have
	\begin{align*}
	M(0)&=\int_{\mathbb{R}}\left|\mathcal{O}_{\alpha}[f](s)\right|^{2} ds-\frac{1}{4}(1+\vert z\vert^2)\int_{\mathbb{R}}|f(t)|^{2}dt=0.
	\end{align*}
	Differentiating \eqref{d} yields
	\begin{align*}
	M^{\prime}(\lambda)=&-|\csc\alpha|^{-\lambda} \ln |\csc\alpha| \int_{\mathbb{R}}|s|^{-\lambda}\left|\mathcal{O}_{\alpha}[f](s)\right|^{2}ds
	-|\csc\alpha|^{-\lambda} \int_{\mathbb{R}}|s|^{-\lambda} \ln |s|\left|\mathcal{O}_{\alpha}[f](s)\right|^{2}ds\\
	&-\frac{1}{4}(1+\vert z\vert^2)C_{\lambda}\int_{\mathbb{R}}|t|^{\lambda} \ln |t||f(t)|^{2}dt-\frac{1}{4}(1+\vert z\vert^2)C_{\lambda}^{\prime}\int_{\mathbb{R}}|t|^{\lambda}|f(t)|^{2}dt.
	\end{align*}
	where constant $C_{\lambda}^{\prime}$ is defined by
	\begin{align*}
	C_{\lambda}^{\prime}= & -\frac{\pi^{\lambda}}{2}\left\{\frac{\Gamma^{2}\left(\frac{1+\lambda}{4}\right) \Gamma\left(\frac{1-\lambda}{4}\right) \Gamma^{\prime}\left(\frac{1-\lambda}{4}\right)+\Gamma^{2}\left(\frac{1-\lambda}{4}\right) \Gamma\left(\frac{1+\lambda}{4}\right) \Gamma^{\prime}\left(\frac{1+\lambda}{4}\right)}{\Gamma^{2}\left(\frac{1+\lambda}{4}\right)}\right\} +\pi^{\lambda} \ln \pi\left\{ \frac{\Gamma^{2}\left(\frac{1-\lambda}{4}\right) }{\Gamma^{2}\left(\frac{1+\lambda}{4}\right)}  \right\} .
	\end{align*}
	Plugging $\lambda=0$ into relation above,  we get
	$C_{0}^{\prime}=\left(\ln \pi-\frac{\Gamma^{\prime}(1 / 4)}{\Gamma(1 / 4)}\right).$
	By virtue of \eqref{Pitt2},  we infer $M(\lambda) \leq 0,$  $\forall \lambda \in[0,1)$ and $M(0)=0$. Therefore, for any $h>0$, we have $M^{\prime}(0+h) \leq 0$, if $h \rightarrow 0$. Hence,
	\begin{align*}
	&- \ln |\csc\alpha| \int_{\mathbb{R}}\left|\mathcal{O}_{\alpha}[f](s)\right|^{2}ds- \int_{\mathbb{R}} \ln |s|\left|\mathcal{O}_{\alpha}[f](s)\right|^{2}ds\\
	&-\frac{1}{4}(1+\vert z\vert^2)C_{0}\int_{\mathbb{R}} \ln |t||f(t)|^{2}dt-\frac{1}{4}(1+\vert z\vert^2)C_{0}^{\prime}\int_{\mathbb{R}}|f(t)|^{2}dt \leq 0.
	\end{align*} 
\end{proof}
\subsection{Local uncertainty principles for $\mathcal{O}_{\alpha}$-tranform}
It is well-known that the classical Heisenberg's uncertainty principle in the Fourier domain points out that if a signal $f$ is well concentrated in the natural domain, then it is impossible for the corresponding Fourier transform $F[f]$ to be well localized around a point in the spectral domain. However, it doesn't preclude $f$ from being localized within a small neighborhood of two or more widely separated points. In fact, the latter phenomenon cannot occur either, and it is the motive of local uncertainty inequalities to make this precise. 
In this part,  we obtain some local uncertainty principles for the $\mathcal{O}_{\alpha}$-transform as follows.
\begin{theorem}
	Let $f \in L^{2}(\mathbb{R})$ and $\lambda \in (0,\frac{1}{2})$. Then, for any $E \subset \mathbb{R}$ with finite measure, the following inequality holds:
\begin{equation*}
	\int_{\mathbb{R}}|t|^{2 \lambda}|f(t)|^{2} dt \geq \frac{4}{(1+\vert z\vert^2)K_{\lambda}|E|^{2 \lambda}} \int_{E}\left|\mathcal{O}_{\lambda}[f](s)\right|^{2} ds.
	\end{equation*}
	where $K_{\lambda}$ is a constant.
\end{theorem}
\begin{proof}
	For a given set $E \subset \mathbb{R}$ with finite measure,  the local uncertainty inequality in the Fourier domain gives us
	\begin{equation*}
	\int_{E}|F[f](s)|^{2} ds \leq K_{\lambda}|E|^{2 \lambda}\left\||t|^{\lambda} f(t)\right\|_{2}^{2}.
	\end{equation*}
	Therefore, 
	\begin{equation*}
	\int_{E}|F[\widetilde{f}](s)|^{2} d s \leq K_{\lambda}|E|^{2 \lambda}\left\||t|^{\lambda} \widetilde{f}(t)\right\|_{2}^{2} .
	\end{equation*}
	Equivalently,
	\begin{equation}\label{Local}
	|\csc\alpha|\int_{E}|F[\widetilde{f}](s\csc\alpha)|^{2} ds \leq K_{\lambda}|E|^{2 \lambda}\left\||t|^{\lambda}\widetilde{f}(t)\right\|_{2}^{2}.
	\end{equation}
	Using the identities $|\widetilde{f}(t)|=|f(t)|$ and $\left|\mathcal{F}_{\alpha}[f](s)\right|=\sqrt{|\csc\alpha|}\left|F[\widetilde{f}](s\csc\alpha)\right|$ in \eqref{Local}, we get
	$$
	\int_{\mathbb{R}}|t|^{2 \lambda}|f(t)|^{2} dt \geq \frac{1}{K_{\lambda}|E|^{2 \lambda}} \int_{E}\left|\mathcal{F}_{\alpha}[f](s)\right|^{2} ds.
	$$
	Similarly,  we infer
	$$
	\int_{\mathbb{R}}|t|^{2 \lambda}|f(t)|^{2} d t \geq \frac{1}{K_{\lambda}|E|^{2 \lambda}} \int_{E}\left|\mathcal{F}_{\alpha}[f](-s)\right|^{2} ds.
	$$
Using \eqref{complex} and $\mathcal{O}_{\alpha}[f](s)=\frac{1}{2}\left[\mathcal{F}_{\alpha}[f](s)+z\mathcal{F}_{\alpha}[f](-s)\right]$,  we infer that
	$$
	\displaystyle\int_{\mathbb{R}}|t|^{2 \lambda}|f(t)|^{2} d t \geq \frac{4}{(1+\vert z\vert^2)K_{\lambda}|E|^{2 \lambda}} \displaystyle\int_{E}\left|\mathcal{O}_{\alpha}[f](s)\right|^{2} ds,
	$$
	which is the desired result.
\end{proof}
\subsection{Hardy's and Beurling-H{\"o}rmander's uncertainty principles}
Towards the end of this section, we formulate another local uncertainty inequality under the Sobolev-type uncertainty inequality. Hardy proposed the uncertainty principle \cite{hardy1933theorem} concerning the decay of a function $f$ and $F[f]$ at infinity.   The Hardy uncertainty principle says that no function is better
localized together with its Fourier transform than the Gaussian, specifically as follows:\\
Let $f \in L^{2}(\mathbb{R})$. If $|f(t)|=O(e^{-\pi \xi t^{2}})$ and $|F[f](s)|=O(e^{-s^{2} / 4 \pi \xi})$ for $\xi>0$,  then $f$ must be of the form
$$
f(t)=K \mathrm{e}^{-\pi \xi t^{2}}, K \in \mathbb{C} 
$$
Now, we formulate  Hardy's uncertainty principle concerning the decay of a function $f$ and the corresponding $\mathcal{O}_{\alpha}$-transform as follows:
\begin{theorem}[Hardy's uncertainty principle for $\mathcal{O}_{\alpha}$-transform]
	Let $f$  be a square-integrable function over $\mathbb{R}$ such that
$	|f(t)|=O(e^{-\pi \xi t^{2}})$ and $\left|\mathcal{O}_{\alpha}[f](s)\right|=O\big(e^{-s^{2} / 4 \pi \xi\sin^2\alpha}\big)$
	for some $\xi>0$. Then
	\begin{equation*}
	f(t)=K \mathrm{e}^{-\left(\pi \xi+i\frac{\cot\alpha}{2} t^{2}\right)}, K \in \mathbb{C}.
	\end{equation*}
\end{theorem}
\begin{proof}
	We remember that $|\widetilde{f}(t)|=|f(t)|$, $\left|\mathcal{F}_{\alpha}[f](s)\right|=\sqrt{|\csc\alpha|} \left|F[\widetilde{f}](s\csc\alpha)\right|$ and 
	\begin{equation}\label{R2}
	\left|\mathcal{O}_{\alpha}[f](\omega)\right|\leq\frac{1}{2}\left[\left|\mathcal{F}_{\alpha}[f](s)\right|+|z|\left|\mathcal{F}_{\alpha}[f](-s)\right|\right].
	\end{equation}
	 Also, note that
	$
	|\widetilde{f}(t)|=O\big(e^{-\pi \xi t^{2}}\big) $ then
	$\left|F[\widetilde{f}](s\csc\alpha)\right|=O\big(e^{-(s\csc\alpha)^{2}/ 4 \pi \xi}\big).
	$
	Therefore	
	 $$\left|\mathcal{O}_{\alpha}[f](\omega)\right|=O\big(e^{-s^{2} / 4 \pi \xi\sin^2\alpha}\big).$$ Using the classical Hardy's principle for the function $\widetilde{f}$,  we receive 
	$\widetilde{f}(t)=K \mathrm{e}^{-\pi \xi t^{2}}, K \in \mathbb{C}.$
	Thus, we obtain
	$
	f(t)=K \mathrm{e}^{-\left(\pi \xi+i\frac{\cot\alpha}{2} t^{2}\right)},$ with $K \in \mathbb{C}.
	$
	The proof is complete.
\end{proof}
A variant of the classical Hardy's uncertainty principle is Beurling's uncertainty principle \cite{beurling1989collected}. A well-known result in \cite{Hormander1991uniqueness}  is the Beurling-Hörmander theorem, which states that: Let $f $ be a 
integrable function on $\mathbb{R}$ with respect to the Lebesgue measure satisfing $F[f] \in L^{1}(\mathbb{R})$ and if
\begin{equation}\label{Beurling_0}
\int_{\mathbb{R}} \int_{\mathbb{R}}|f(t)||F[f](\omega)| e^{|ts|} dt  d s<+\infty,
\end{equation}
then $f=0$ almost everywhere. A strong multidimensional version of this theorem has been established by Bonami, Demange
and Jaming in \cite{Bonami2003hermite}, who have shown that if $f$ is a square-integrable function on $\mathbb{R}^n$ with respect to the Lebesgue measure, then $$\int_{\mathbb{R}^n}\int_{\mathbb{R}^n}\frac{|f(x)|\ |\widetilde{f}(y)|}{(1+|x|+|y|)^d}e^{|\langle x\textfractionsolidus y\rangle|}dxdy<+\infty$$ with $d\geq 0,$ if and only if f can be written as $f(x)=P(x) e^{-\langle Ax\textfractionsolidus x\rangle}$, where $A$ is a real positive definite symmetric matrix and $P$ is a polynomial with  $\operatorname{deg} (P)<\frac{d-n}{2}.$ In particular for $d\leq n$, then $f$ is identically zero.
We will give a variant version of theorem in \cite{Hormander1991uniqueness} for the $\mathcal{O}_{\alpha}$ transformation as follows:
\begin{theorem}[Beurling-H{\"o}rmander's theorem for $\mathcal{O}_{\alpha}$-transform]
	Assume that  $f$ be a integrable function on $\mathbb{R}$ with respect to the Lebesgue measure satisfing $\mathcal{O}_{\alpha}[f] \in L^{1}(\mathbb{R})$ and
	\begin{equation}\label{Beurling_1}
	\int_{\mathbb{R}} \int_{\mathbb{R}}\left|f(t)\right|\left|\mathcal{O}_{\alpha}[f](\omega)\right| e^{|t \omega / b|}  dt ds<+\infty,
	\end{equation}
	then $f=0$ almost everywhere.
\end{theorem}
\begin{proof}
With the help of relation $$|\widetilde{f}(t)|=|f(t)|\quad  \text{and}\quad  \left|\mathcal{F}_{\alpha}[f](s)\right|=\sqrt{|\csc\alpha|} \left|F[\widetilde{f}](s\csc\alpha)\right|$$ and utilizing inequality \eqref{Beurling_0},  we deduce that
	\begin{equation*}
	\int_{\mathbb{R}} \int_{\mathbb{R}}|\widetilde{f}(t)||\mathcal{F}_{\alpha}[f](s)| e^{|t \omega|} dt  ds=\sqrt{|\csc\alpha|} \int_{\mathbb{R}} \int_{\mathbb{R}}|\widetilde{f}(t)|\left|F[\widetilde{f}](s\csc\alpha)\right| e^{|t \omega / b|}  dt  ds<+\infty.
	\end{equation*}
	Similarly
	\begin{equation*}
	\int_{\mathbb{R}} \int_{\mathbb{R}}|\widetilde{f}(t)||\mathcal{F}_{\alpha}[f](-s)| e^{|t \omega|} dt  ds<+\infty.
	\end{equation*}	
	Thus,  by using \eqref{R2}, we directly obtain \eqref{Beurling_1}.  Moreover, we then have $\widetilde{f}(t)=0$.  Hence $f$ is identically zero.
\end{proof}

\vskip 0.3cm
\noindent \textbf{Acknowledgements}\\
The author expresses sincere gratitude to the anonymous referees who provided many helpful suggestions for this manuscript.\\
\noindent \textbf{Disclosure statement}\\
\noindent No potential conflict of interest was reported by the author.\\
\noindent \textbf{Ethics declarations}\\
\noindent We confirm that all the research meets ethical guidelines and adheres to the legal requirements of the study country.\\
\noindent \textbf{Funding}\\
\noindent This research received no specific grant from any funding agency.\\
\noindent \textbf{ORCID}\\
\noindent \textsc{Trinh Tuan}  \url{https://orcid.org/0000-0002-0376-0238}\\
\noindent \textsc{Lai Tien Minh} \url{https://orcid.org/0000-0003-2656-8246}\\

\end{document}